\documentstyle[amstex,amscd]{article}
\title{Projectivity criterion of Moishezon spaces 
and density of projective symplectic varieties}
\author{Yoshinori Namikawa}
\date{ }
\begin{document}
\maketitle

\begin{center}
{\bf Abstract}. 
\end{center}

A Moishezon manifold is a projective manifold 
if and only if it is a K{\"a}hler manifold [Mo 1]. 
However, a singular Moishezon space is not 
generally projective even if it is a K{\"a}hler 
space [Mo 2]. Vuono [V] has given a 
projectivity criterion for Moishezon spaces 
with isolated singularities. In this paper 
we shall prove that a Moishezon space with 
1-rational singularities is projective 
when it is a K{\"a}hler space (Theorem 6). 

We shall use Theorem 6 to show the 
density of projective symplectic varieties 
in the Kuranishi family of a 
(singular) symplectic variety 
(Theorem 9), which is a generalization of 
the result by Fujiki [Fu 1, Theorem 4.8] 
to the singular case.

In the Appendix we give a supplement 
and a correction to the 
previous paper [Na] where singular symplectic 
varieties are dealt with.  
\vspace{0.12cm}

{\em Acknowledgement}. The author would like to 
thank K. Oguiso, who asked him if Theorem 9 
holds in the singular case.  
\vspace{0.15cm}

\begin{center}
{\bf 1. Projectivity of a K{\"a}hler 
Moishezon variety }
\end{center}

In this paper, a compact complex variety 
means a compact, irreducible and reduced 
complex space. We mean by a Moishezon variety $X$ 
a compact complex variety $X$ with $n$ 
algebraically independent meromorphic functions, 
where $n = \dim X$, [Mo 1]. 

For a complex space $X$ we say that $X$ is 
K{\"a}hler 
if $X$ admits a K{\"a}hler metric (form) 
in the sense 
of [Gr], [Mo 2] (cf.[B]). \vspace{0.2cm}

{\bf Definition 1}. Let $X$ be a compact complex 
variety. 

(1) An element $b \in H_{2}(X, {\bold Q})$ 
is an {\em analytic homology class} 
if $b$ is represented by a 2-cycle 
$\Sigma \alpha_j C_j$ where $C_j$ are 
complex subvarieties of dimension $1$ and    
$\alpha_j \in {\bold Q}$. Denote 
by $A_{2}(X, {\bold Q})$ the 
subspace of $H_2(X, {\bold Q})$ spanned 
by analytic homology classes. 
Define $A_2(X, {\bold R}) := A_2(X, {\bold Q})
\otimes_{\bold Q}{\bold R}$. \vspace{0.12cm}
 
(2) For an element $b \in A_2(X, {\bold Q})$ 
and a line bundle $L$ of $X$, the intersection 
number $(b, L)$ is well-defined. Two elements 
$b$, $b'$ of $A_2(X, {\bold Q})$ are said to 
be numerically equivalent if 
$(b, L) = (b', L)$ for all line bundles $L$ 
on $X$. Denote by $N_1(X)_{\bold Q}$ the 
quotient ${\bold Q}$-vector space of 
$A_2(X, {\bold Q})$ by this numerical 
equivalence. Define $N_1(X)_{\bold R} 
:= N_1(X)_{\bold Q}\otimes{\bold R}$. 
\vspace{0.12cm}

(3) Define $\mathrm{Pic}(X)_{\bold Q} 
:= \mathrm{Pic}(X)\otimes_{\bold Z}{\bold Q}$
and $\mathrm{Pic}(X)_{\bold R} 
:= \mathrm{Pic}(X)\otimes_{\bold Z}{\bold R}$.
Two line bundles $L$ and $L'$ are 
numerically equivalent if $(b, L) 
= (b, L')$ for all $b \in A_2(X, {\bold Q})$. 
Let $N^1(X)$ be the abelian 
group of numerical classes of 
line bundles on $X$. 
Define $N^1(X)_{\bold Q} := N^1(X)\otimes
_{\bold Z}{\bold Q}$ and 
$N^1(X)_{\bold R} := N^1(X)\otimes_{\bold 
Z}{\bold R}$.
   \vspace{0.15cm}

{\bf Proposition 2}([Ko-Mo, (12.1.5)]) 
{\em 
Let $X$ be a Moishezon variety 
with 1-rational singularities, that 
is, $X$ is normal and has a 
resolution $\pi: Y \to X$ 
such that $R^1\pi_*{\cal O}_Y = 0$. 
Then an analytic homology 
class $b \in A_2(X, {\bold Q})$ is 
zero if it is numerically equivalent to $0$. 
In particular, 
$A_2(X, {\bold Q}) = N_1(X)_{\bold Q}$}. 
\vspace{0.15cm}

{\em Sketch of proof}: Let $\pi : Y \to X$ 
be a resolution such that $Y$ is projective. 
Let $N_1(Y/X)_{\bold Q}$ be the subspace 
of $N_1(Y)_{\bold Q}$ generated by the 
classes of curves contained in a fiber of 
$\pi$. We have an exact sequence 
$$  0 \to N_1(Y/X)_{\bold Q} 
\to N_1(Y)_{\bold Q} \to N_1(X)_{\bold Q} 
\to 0.$$
We use the condition $R^1\pi_*{\cal O}_Y 
= 0$ to prove the middle exactness of the 
sequence. 

We may assume that $b$ is represented by a  
curve $C$ on $X$. By cutting out $\pi^{-1}(C)$ 
by general hyperplane sections of $Y$, we can find 
a curve $C'$ such that, as a cycle, 
$\pi_*(C') = mC$ for 
a positive integer $m$. Define $b' := [(1/m)C'] 
\in N_1(Y)_{\bold Q}$. By the exact sequence, 
$b' \in N_1(Y/X)_{\bold Q}$. 
Thus $b'$ is represented 
by a ({\bf Q})curve $D$ contained in 
some fibers of 
$\pi$. By definition, $D$ and $(1/m)C'$ are 
numerically equivalent. 
Since numerical equivalence 
and homological equivalence coincide on $Y$ 
(cf. [Mo 1, page 83, Theorem 9]), 
these are also homological equivalent. 
Therefore, $\pi_*(D)$ and $\pi_*((1/m)C')$ are 
homologically equivalent on $X$. 
Since $[\pi_*(D)] = 0$, $b = [\pi_*((1/m)C')] 
= 0$. \vspace{0.15cm}       

{\bf Lemma 3}. {\em Let $X$ be a K{\"a}hler, 
Moishezon 
variety with a K{\"a}hler form $\omega$. 
Assume that 
an element $L \in \mathrm{Pic}(X)_{\bold Q}$ 
satisfies the 
following condition $(**)$}. 

(**) {\em For every curve $C \subset X$}, 
$$  \int_{C}\omega = (C. L).  $$ 

{\em Then, for any subvariety $W \subset X$ of 
dimension $k$, ${(L^k)_W} > 0$}.  
\vspace{0.15cm}

{\em Proof}. By the assumtion, 
$L$ is nef, hence 
${(L^k)_W} \geq 0$ (cf. [Kl]). We shall derive 
a contradiction by assuming that 
${(L^k)_W} = 0$. 

Put $\omega' := \omega\vert_W$. 
Then $\omega'$ becomes 
a K{\"a}hler form of $W$, and $W$ is a 
K{\"a}hler Moishezon variety. 
Let $p \in W$ be a smooth 
point, and let $h: \hat W \to W$ 
be the blowing up 
at $p$. Denote by $E$ the inverse 
image $h^{-1}(p)$. 
$E$ is a Cartier divisor of $\hat W$. 
By using a 
Hermitian metric on ${\cal O}_{\hat W}(E)$, 
one can 
define a $d$-closed $(1,1)$-form $\alpha_E$ 
on $\hat W$ in such a way that 
$\omega_{\epsilon} := 
h^*\omega - \epsilon \alpha_E$ 
becomes a K{\"a}hler form 
on $\hat W$ if $\epsilon >0$ is 
sufficiently small. 

Fix a sufficiently small rational number 
$\epsilon > 0$. Put $L' := L\vert_W$ and 
$F:= h^*{L'}-{\epsilon}E$. 
By the condition 
$(**)$, for every curve $C \subset \hat W$, 
we have $(C. F) = \int_{C}
(h^*{\omega'} - \epsilon \alpha_E)$.  

In particular, $F$ is a nef {\bf Q}-line bundle on 
$\hat W$. 
Therefore, $(F)^k \geq 0$. 

On the other hand, since 
$(h^*{L'})^i(E)^{k-i} = 0$ 
for $i \ne 0, k$, we have 
$(F)^k = (L')^k + \epsilon^k
(-E)^k$. By the assumption, $(L')^k = 0$. 
By the 
definition of $E$, 
we have $(-E)^k < 0$. Hence 
$(F)^k < 0$, a contradiction. \vspace{0.2cm}

{\bf Theorem 4} 
(cf. [Mo 1, page 77, Theorem 6]) 
{\em Let $X$ be a Moishezon variety. 
Assume that 
an element $L \in \mathrm{Pic}(X)_{\bold Q}$ 
satisfies the 
inequality $(L)^k_W > 0$ for every 
$k$ dimensional 
subvariety $W$ of $X$. Then $L$ is ample.} 
\vspace{0.2cm}

By Lemma 3 and Theorem 4 we have the 
following corollary. \vspace{0.15cm}

{\bf Corollary 5}. {\em Let $X$ be a 
K{\"a}hler Moishezon variety 
with a K{\"a}hler 
form $\omega$. Assume that an element 
$L \in \mathrm{Pic}(X)_{\bold Q}$  
satisfies the equality for 
any curve} $C \subset X$:  

$$ (C. L) = \int_{C}\omega. $$

{\em Then $L$ is ample.} \vspace{0.2cm}

{\bf Theorem 6}. {\em Let $X$ be a Moishezon 
variety with 1-rational singularities (cf. 
Proposition 2). 
If $X$ is K{\"a}hler, then $X$ is projective.}  
\vspace{0.15cm}

{\em Proof}. 

Since the numerical equivalence and the 
homological equivalence coincide for 
(analytic) 1-cycle by Proposition 2, 
we have   
a natural map $\alpha : N^1(X)_{\bold Q}  
\to (A_2(X, {\bold Q}))^*$ and $\alpha$ 
is an isomorphism.  

Taking the tensor product with 
${\bold R}$, we have  
a map $\alpha_{\bold R} : 
N^1(X)_{\bold R} 
\to (A_2(X, {\bold R}))^*$ and 
$\alpha_{\bold R}$ is an isomorphism. 
By the 2-nd cohomolgy class defined 
by the K{\"a}hler form $\omega$ 
(cf. [B, (4.15)]) 
one can regard the K{\"a}hler 
form as an element of 
$(A_2(X, {\bold R}))^*$. 
Since $\alpha_{\bold R}$ is surjective, 
there is an element   
$d \in N^1(X)_{\bold R}$  
such that $(C.d) = \int_{C}\omega$ 
for every curve $C \subset X$. 

Approximate $d \in N^1(X)_{\bold R}$ 
by a convergent sequence $\{d_m\}$ of 
rational elements 
$d_m \in N^1(X)_{\bold Q}$.  
 
Let us fix the basis $b_1, ..., b_l$ 
of the vector space 
$N^1(X)_{\bold Q}$. 
Each $b_i$ is represented by an 
element   
$B_i \in \mathrm{Pic}(X)$. 
Now $d$ (resp. $d_m$) is 
represented by an element in $\mathrm{Pic}
(X)_{\bold R}$ 
(resp. $\mathrm{Pic}(X)_{\bold Q}$) 
$D:= \Sigma x_i B_i$ (resp. $D_m := 
\Sigma x^{(m)}_i B_i$) 
such that $\lim x^{(m)}_i = x_i$. 

Put $E_m := D_m - D$. Then there are 
$d$ closed $(1,1)$-forms $\alpha_m$ 
corresponding to $E_m$ such that  
$\{\alpha_m\}$ uniformly converge to $0$. 
     
If $m$ is chosen sufficiently large, 
then $\omega_m := \omega + \alpha_m$ 
is a K{\"a}hler form. Since 

$$ (C. D_m) = \int_{C}\omega_m $$ 
for every curve $C \subset X$, 
we see that $D_m$ is ample by 
Corollary 5. \vspace{0.15cm}

{\bf Corollary 6'}. {\em Let $X$ be a Moishezon 
variety with rational singularities. 
If $X$ is K{\"a}hler, then $X$ is projective.}
\vspace{0.15cm}

{\bf Remark}. If we do not assume 
that $X$ has 1-rational singularities, 
Theorem 6 is no longer 
true (cf. [Mo 2]). \vspace{0.15cm}

\begin{center}
{\bf 2. Application: Density of projective 
symplectic varieties} 
\end{center} 

A symplectic variety is a compact 
normal K{\"a}hler space $X$  
with the following properties: 
(1) The regular part $U$ 
of $X$ has an everywhere non-degenerate 
holomorphic 
2-form $\Omega$, and 
(2) for a (any) resolution of 
singularities $f : {\tilde X} \to X$ 
such that 
$f^{-1}(U) \cong U$, 
the 2-form $\Omega$ extends to 
a holomorphic 2-form on ${\tilde X}$.  
Here the extended 
2-form may possibly degenerate 
along the exceptional locus.
By definition, $X$ has only 
canonical singularities, 
hence has only rational singularities. 

If $X$ has a resolution 
$f: {\tilde X} \to X$ such that 
$\Omega$ extends to an 
everywhere non-degenerate 2-form 
on ${\tilde X}$, then we say that 
$X$ has a symplectic 
resolution. 

Symplectic varieties with no symplectic 
resolutions are constructed as symplectic 
V-manifolds in [Fu 1]. 
Recently, O'Grady [O] has 
constructed such varieties as the moduli spaces 
of semi-stable torsion free sheaves on a 
polarized K3 surface (cf. [Na, Introduction]). 
His examples are no more V-manifolds. 

These examples satisfy the following 
condition: 

(*): The natural restriction map 

$$  H^2(X, {\bold Q}) \cong H^2(U, {\bold Q}) $$ 
is an isomorphism.  

In [Na] we have formulated the local 
Torelli problem for these symplectic 
varieties, and proved it. More precisely,
we have proved it for a symplectic variety 
$X$ with the following properties.
 
(a): $\mathrm{Codim}(\Sigma \subset X) 
\geq 4$, where $\Sigma := \mathrm{Sing}(X)$, 
\vspace{0.12cm}
 
(b): $h^1(X, {\cal O}_X) = 0$,  
$h^0(U, \Omega^2_U) = 1$, and \vspace{0.12cm}

(c): (*) is satisfied.\footnote{
Note that this condition is equivalent to the 
condition $(*)$ in [Na, Remark (2)]
(cf. (b) in the proof of [Na, Prop. 9]).}  
    \vspace{0.2cm}

Let $X$ be a symplectic variety satisfying 
(a), (b) and (c). Let $0 \in S$ be the Kuranishi 
space of $X$ and $\bar{\pi}: {\cal X} 
\to S$ be the universal family such that 
$\bar{\pi}^{-1}(0) = X$.  Let ${\cal U}$ 
be the locus in ${\cal X}$ where 
$\bar{\pi}$ is a smooth map. We denote 
by $\pi$ the restriction $\bar{\pi}$ 
to ${\cal U}$. $S$ is nonsingular 
by the condition (a). 
Note that every fiber of $\bar{\pi}$ is 
a symplectic variety satisfying 
(a), (b) and (c) (cf. [Na]). 

The cohomology $H^2(U, {\bold C})$ 
admits a natural mixed Hodge structure 
because $U$ is a Zariski open 
subset of a compact K{\"a}hler space 
(cf. [Fu 2]). 
In our case, because of condition (a), 
it is pure of weight 2, and 
the Hodge decomposition is given by 

$$ H^2(U, {\bold C}) = H^0(U, \Omega^2_U) 
\oplus H^1(U, \Omega^1_U) \oplus 
H^2(U, {\cal O}_U). $$

For details of this, see the footnote 
of the e-print version of [Na, p.21]. 
Moreover, the tangent space $T_{S,0}$ 
at the origin is canonically isomorphic 
to $H^1(U, \Theta_U)$. By a holomorphic 
symplectic 2-form $\Omega$, $H^1(U, \Theta_U)$ 
is identified with $H^1(U, \Omega^1_U)$. 

Let us fix a resolution $\tilde X$ of $X$. 
For $\alpha \in H^2(U, {\bold C})$, denote 
by $\tilde\alpha \in H^2(\tilde X, {\bold C})$ 
the image of $\alpha$ by the map 
$H^2(U, {\bold C}) 
\cong H^2(X, {\bold C}) \to H^2(\tilde X, 
{\bold C})$. The holomorphic symplectic 
2-form $\Omega$ defines a holomorphic 2-form 
on $\tilde X$ by the definition of a symplectic 
variety. We denote by the same symbol $\Omega$ 
this holomorphic 2-form on ${\tilde X}$. 
Here we normalize $\Omega$ so that 
$\int_{\tilde X}(\Omega {\overline \Omega})^l 
= 1$.
 
One can define 
a quadratic form $q$ on $H^2(U, {\bold C})$ by 

$$q(\alpha) :=  
l/2 {\int_{\tilde X}
(\Omega {\overline \Omega})^{l-1}
\tilde\alpha^2} 
+ (1-l)\int_{\tilde X}\Omega^l
{\overline\Omega}^{l-1}
\tilde\alpha
\int_{\tilde X}\Omega^{l-1}
\overline\Omega^l\tilde\alpha, $$ 
where $\dim X = 2l$. 
Note that $q$ is independent of the choice 
of the resolution $\tilde X$. 
The quadratic form  
$q$ is defined over $H^2(U, {\bold R})$. 
By the same argument as [Be, 
Th{\'e}or{\`e}me 5, (a), (c)], 
we can write $q(\alpha) 
= c_1 ({\int_{\tilde X}[f^*\omega]^{2l-2}
{\tilde
\alpha}^2}) - c_2 
({\int_{\tilde X}[f^*\omega]^{2l-1}
{\tilde\alpha}})^2$ by a K{\"a}hler 
form $\omega$ 
on $X$ and by suitable positive 
real constants $c_1$ and $c_2$. 
Let us define 
$H^2_0(U, {\bold R}) := \{\alpha \in 
H^2(U, {\bold R}); 
{\int_{\tilde X}[f^*\omega]^{2l-1}
{\tilde \alpha}} = 0\}$. 
Then we have a direct 
sum decomposition $H^2(U, {\bold R}) 
= 
H^2_0(U, {\bold R})\oplus{\bold R}[\omega]$, 
where 
$H^2_0(U, {\bold R})$ and ${\bold R}[\omega]$ 
are orthogonal with respect to $q$. 

We shall 
prove that 
$q([\omega]) > 0.$ 

It is easily 
checked that 
$q([\omega]) = l/2 {\int_{\tilde X}
(\Omega {\overline \Omega})
^{l-1}[f^*\omega]^2}$. 
We assume that the resolution 
$f: \tilde X \to X$ 
is obtained by a succession of blowing ups with 
smooth centers contained in the singular locus. 
Let $\{E_i\}$ be the exceptional divisors of 
$f$. Then, for sufficiently 
small positive real numbers $\epsilon_i$, 
$[f^*\omega] - \Sigma \epsilon_i [E_i]$ 
is a K{\"a}hler class on $\tilde X$ 
(cf. [Fu 3, Lemma 2]). 
By [W, Corollaire au 
Th{\'e}or{\`e}me 7, p.77] 
we have $l/2 {\int_{\tilde X}
(\Omega 
{\overline \Omega})^{l-1}([f^*\omega] - 
\Sigma \epsilon_i [E_i])^2} 
> 0$.   
On the other hand, by [Na, Remark (1), p. 24] 
we have ${\int_{\tilde X}
(\Omega 
{\overline \Omega})^{l-1}[f^*\omega][E_i]} = 
{\int_{\tilde X}
(\Omega {\overline \Omega})^{l-1}[E_i]^2} = 0$; 
hence we have $q([\omega]) > 0$.      

Denote by $Q : 
H^2(U, {\bold C}) \times H^2(U, {\bold C}) 
\to {\bold C}$ the symmetric 
bilinear form defined by $q$. 
With respect to $Q$, $H^0(U, \Omega^2_U)
\oplus H^2(U, {\cal O}_U)$ is orthogonal 
to $H^1(U, \Omega^1_U)$. 
Let us define 
$N := \{v \in H^1(U, \Omega^1_U); 
Q(v, x) = 0$ for any 
$x \in H^1(U, \Omega^1_U)\}.$  
Since $q([\omega]) > 0$ 
for a K{\"a}hler form $\omega$ 
on $X$, $N$ does not coincide with 
$H^1(U, \Omega^1_U)$. By the identification 
of $H^1(U, \Omega^1_U)$ with $T_{S,0}$, we 
regard $N$ as a subspace of $T_{S,0}$.     

Later we shall prove that the quadratic 
form $q$ is non-degenerate and, in fact,  
$N = 0$. But, before doing this, we 
first prove the density of projective 
symplectic varieties in a rather imcomplete 
form (cf. Proposition 7 below). 
After that we will show that 
$q$ is non-degenerate by using  
Proposition 7. As a consequence, we 
will see that, in Proposition 7, the 
assumption for $T_{S_1,0}$ is not 
necessary; hence we can prove  
the density in a complete form (cf. 
Theorem 9).  
\vspace{0.15cm}

{\bf Proposition 7}. {\em Notation and 
assumptions being the same as above, 
let $0 \in S_1 \subset S$ be a 
positive dimensional non-singular subvariety 
of $S$ such that $T_{S_1,0}$ is not contained 
in $N$. Then, for any open 
neighborhood $0 \in V \subset S$, 
there is a point $s \in V \cap S_1$ 
such that ${\cal X}_s$ is a projective 
symplectic variety.}  \vspace{0.15cm}

{\em Proof}. We may assume that 
$\dim S_1 = 1$. Denote by ${\cal X}_1$ 
the fiber product ${\cal X}\times_S S_1$ 
and denote by $\bar{\pi}_1$ the induced 
map from ${\cal X}_1$ to $S_1$.  
Take a resolution of singularities 
$\nu: {\cal Y}_1 \to {\cal X}_1$ in 
such a way that $\nu$ is an isomorphism 
over smooth locus of ${\cal X}_1$. 
We also assume that $\nu$ is obtained by 
the succession of blowing ups with smooth centers. 
So there exists a $\nu$-ample 
divisor of the form 
$- \Sigma \epsilon_i {\cal E}_i$, where 
${\cal E}_i$ are $\nu$-exceptional 
divisors and $\epsilon_i$ 
are positive rational numbers.   

Let $S_1^0$ be the set of points 
$s \in S_1$ where 
$(\bar{\pi}_1\circ \nu)^{-1}(s)$ are smooth. 
Then $S_1^0$ is a non-empty Zariski 
open subset of $S_1$. 

We may assume that $0 \in S_1^0$. 
In fact, if $0 \notin S_1^0$, then 
we take a point $s \in S_1^0 \cap V$. 
Then the family $\bar{\pi}: {\cal X} \to S$ 
can be regarded as the Kuranishi 
family of ${\cal X}_s$ near $s \in S$ 
because $H^0(X, \Theta_X) = 0$ by the 
condition (b). For this point $s \in S$, 
${\cal X}_s$ satisfies all conditions 
(a), (b) and (c). So, if the theorem 
holds for ${\cal X}_s$, then we can 
find a point $s' \in S_1^0 \cap V$ 
where ${\cal X}_{s'}$ is projective. 

In the remainder we shall assume that 
$0 \in S_1^0$. Thus, if $S$ is chosen 
sufficiently small, then $\nu : {\cal Y}_1 
\to {\cal X}_1$ is a simultaneous resolution 
of $\{{\cal X}_{1,s}\}$, $s \in S_1$.     
\vspace{0.15cm}
   
{\bf Claim} {\em For any open 
neighborhood $0 \in V \subset S$, 
there is a point $s \in V \cap S_1$ 
such that ${\cal Y}_{1,s}$ is a projective 
variety.} \vspace{0.15cm}

If the claim is justified, then, for 
such a point $s$, ${\cal X}_{1,s}$ 
is a Moishezon variety. On the other hand, 
${\cal X}_{1,s}$ is a symplectic variety 
satisfying (a), (b) 
and (c) (cf. [Na]). In particular, 
${\cal X}_{1,s}$ is a 
K{\"a}hler Moishezon variety 
with rational singularities. 
By Theorem 6, we conclude 
that ${\cal X}_{1,s}$ is projective. 
\vspace{0.15cm}

{\em Proof of Claim}. 

(i): Put $Y ={\cal Y}_{1,0}$. 
By definition of ${\cal X}_1$, ${\cal X}_{1,0} 
= X$. The bimeromorphic map $\nu_0 : Y \to X$ 
is a resolution of singularities. Put $E_i := 
{\cal E}_{i,0}$ where ${\cal E}_i$ are 
$\nu$-exceptional divisors. By the construction 
of $\nu$, there are positive rational 
numbers $\epsilon_i$ such that 
$- \Sigma \epsilon_i E_i$ is $\nu_0$-ample. 
\vspace{0.12cm}
 
(ii):  We have a constant sheaf 
$R^2{\bar{\pi}}_*{\bold C}$ 
on $S$. There is an isomorphism 
$R^2{\bar{\pi}}_*{\bold C}\otimes_{\bold C}
{\cal O}_S \cong 
R^2{\pi}_*{\bold C}\otimes_{\bold C}
{\cal O}_S$.  
The right hand side is filtered 
as $R^2{\pi}_*{\bold C}\otimes_{\bold C}
{\cal O}_S = {\cal F}^0 \supset {\cal F}^1 
\supset {\cal F}^2 \supset 0$ in such a way that 
$Gr^i_{\cal F} =  
 R^{2-i}\pi_*
\Omega^i_{{\cal U}/S}.$   

Over each point $s \in S$, this filtration  
gives the Hodge decomposition of 
$H^2({\cal X}_s)$. The natural 
mixed Hodge structure on $H^2({\cal X}_s)$ 
is pure, and its $(i,j)$-component 
$H^{i,j}({\cal X}_s)$ $(i+j = 2)$ 
is given by  
$H^j({\cal U}_s, \Omega^i_{{\cal U}_s})$.

For $a \in H^2(X, {\bold R}) = 
\Gamma(S, R^2{\bar{\pi}}_*{\bold R})$, 
define $S_a$ to be the locus in $S$ 
where $a \in H^{1,1}({\cal X}_s)$. 

Let $\omega$ be a K{\"a}hler form on $X$ 
and denote by $[\omega] \in 
H^2(X, {\bold R}) (= H^2(U, {\bold R}))$ its 
cohomology class (cf. [B, (4.15)]).  
Then the tangent space $T_{S_{[\omega]},0}$ 
is isomorphic to $\{v \in H^1(U, \Omega^1_U); 
Q(v, [\omega]) = 0 \}$, where we identify 
$T_{S,0}$ with $H^1(U, \Omega^1_U)$ as 
explained above.  Let $v_1 \in 
H^1(U, \Omega^1_U)$ be a generator of 
1-dimensional vector space $T_{S_1,0}$. 

Define a linear map $Q_{v_1}: H^1(U, \Omega^1_U) 
\to {\bold C}$ by $Q_{v_1}(x) = Q(v_1, x)$. 
By the assumption, 
this map is surjective. 
By definition, $\mathrm{Ker}(Q_{v_1})$ 
is a complex hyperplane in $H^1(U, \Omega^1_U)$. 
Therefore, $H^1(U, \Omega^1_U) \cap 
H^2(U, {\bold R})$ is not contained in 
$\mathrm{Ker}(Q_{v_1})$. Hence  
we can find an element $a$ in any 
small open neighborhood of $[\omega] 
\in H^{1,1}(X) \cap H^2(X, {\bold R})$ 
in such a way that $S_a$ intersects 
$S_1$ in $0$ transversely. \vspace{0.12cm}

(iii): Let $\alpha_{E_i}$ be a 
d-closed $(1,1)$ form on $Y$ corresponding to 
$E_i$. If necessary, by multiplying to 
$\{\epsilon_i\}$ a 
sufficiently small positive rational number       
simultaneously, the d-closed (1,1) form 

$$ (\nu_0)^*\omega - \Sigma 
\epsilon_i \alpha_{E_i} $$ 
becomes a K{\"a}hler form on $Y$ (cf. 
[Fu 3, Lemma 2]). 

Take $a \in H^{1,1}(X) \cap H^2(X, {\bold R})$ 
as in the final part of (ii). Then we can find 
a suitable d-closed $(1,1)$-form $\omega'$ on 
$Y$ which represents $(\nu_0)^*a \in 
H^{1,1}(Y) \cap H^2(Y, {\bold R})$ in such   
a way that 

$$ \gamma := \omega' - \Sigma \epsilon_i 
\alpha_{E_i} $$ 
becomes a K{\"a}hler form on $Y$. 
 
$[E_i]$ remains of $(1,1)$ type in 
$H^2({\cal Y}_{1,s})$ for arbitrary $s \in S_1$ 
because the Cartier divisor $E_i$ extends 
sideways in the family ${\cal Y}_1 \to S_1$. 
Hence, by the choice of $a$, $(\nu_0)^*a 
- \Sigma \epsilon_i [E_i]$ is no more 
of type $(1,1)$ for $s \ne 0$ sufficiently near  
$0$. \vspace{0.12cm}

(iv): The last step is the same as 
the proof of [Fu 1, Theorem 4.8, (2)]. 
Approximate $a \in H^{1,1}(X) \cap 
H^2(X, {\bold R})$ by 
a sequence of rational classes 
$\{a_m\}$ ($a_m \in H^2(X, {\bold Q})$).    
   
We put $b := (\nu_0)^*a - 
\Sigma \epsilon_i [E_i]$ 
and $b_m := (\nu_0)^*{a_m} - 
\Sigma \epsilon_i [E_i]$. 

We take a $C^{\infty}$ family of 
K{\"a}hler forms 
$\{\gamma_s\}$ on $\{{\cal Y}_{1,s}\}$ such 
that $\gamma_0 = \gamma$. For each $s \in S_1$, 
the cohomolgy class $b_m$ is represented by a 
unique harmonic 2-form 
$(\omega_m)(s)$ on ${\cal Y}_{1,s}$ 
with respect to $\gamma_s$.  
If $m$ is large, then $(\omega_m)(s)$ becomes 
a K{\"a}hler form on ${\cal Y}_{1,s}$ 
for some point 
$s \in V \cap S_1$. Since $b_m$ is a 
rational cohomolgy class, ${\cal Y}_{1,s}$ 
is a projective manifold by [Ko] for 
this $s$. \vspace{0.15cm}

{\bf Corollary 8}. {\em Let $X$ be a symplectic 
variety of dim $n = 2l$ satisfying (a), (b) and (c). 
Let $q$ 
be the quadratic form on $H^2(U, {\bold R})$ 
defined above. Then $q$ is non-degenerate 
and has signature $(3, B-3)$ where $B := 
\dim H^2(U, {\bold R})$.} \vspace{0.12cm}

{\em Proof}. We first prove this collorary 
when $X$ is projective. Let $\omega$ be a 
K{\"a}hler form that comes from a very 
ample line bundle $L$ on $X$. Take general 
global sections $t_1$, ..., $t_{n-2}$ of 
$L$. We put $T_0 = X$ and, for $1 \leq i \leq n-2$, 
define $T_i \subset X$ to 
be the common 
zeros of $t_1$, ..., $t_i$. Put $\Sigma_i 
= \Sigma \cap T_i$, where $\Sigma := 
\mathrm{Sing}(X)$. By the condition (a), 
$\Sigma_{n-2} = \emptyset$ and 
$T_{n-2}$ is a nonsingular surface. 
For simplicity we put $T := T_{n-2}$.
By [H, Theorem 2], the pair $(T_i - \Sigma_i, 
T_{i+1} - \Sigma_{i+1})$ is $n-i-1$-connected. 
In particular, the restriction maps 
$H^2(T_i\setminus\Sigma_i, {\bold R}) \to 
H^2(T_{i+1}\setminus\Sigma_{i+1}, {\bold R})$ 
are all injective. Hence $H^2(U, 
{\bold R}) \to H^2(T, {\bold R})$ 
is an injection. Note that both sides 
have pure Hodge structures of weight 
2 (cf. the footnote on p.21 of the e-print 
version of [Na]) and this injection 
is a morphism of Hodge structures.  
Define $H^2_0(T, {\bold R}) 
:= \{\alpha \in H^2(T, {\bold R}); 
(\omega\vert_T).\alpha = 0 \}$. Then 
we have the restriction 
map $H^2_0(U, {\bold R}) \to H^2_0(T, {\bold R})$, 
which is an injection. This restriction map 
induces 
an injection $H^{1,1}_0(U)_{\bold R} \to 
H^{1,1}_0(T)_{\bold R}$, where 
$H^{1,1}_0(U)_{\bold R} := H^{1,1}(U) \cap 
H^2_0(U, {\bold R})$ and $H^{1,1}_0(T)_{\bold R} 
:= H^{1,1}(T) \cap H^2_0(T, {\bold R})$. 
 Let $q'$ be the quadratic 
form on $H^{1,1}_0(T)_{\bold R}$ defined by the cup 
product. By the definition of $H^2_0(U, {\bold R})$, 
$q\vert_{H^{1,1}_0(U)_{\bold R}} = 
c_1{q'}\vert_{H^{1,1}_0(U)_{\bold R}}$ for 
a suitable positive real constant $c_1$. 
Since ${q'}\vert_{H^{1,1}_0(U)_{\bold R}}$ 
is negative-definite, 
$q\vert_{H^{1,1}_0(U)_{\bold R}}$ is 
also negative-definite. We have a direct 
sum decomposition with respect to $q$: 
$$ H^2(U, {\bold R}) = {\bold R}[\omega] 
\oplus (H^{2,0}(U)\oplus H^{0,2}(U))_{\bold R} 
\oplus H^{1,1}_0(U)_{\bold R}, $$ 
where $(H^{2,0}(U)\oplus H^{0,2}(U))_{\bold R} 
:= (H^{2,0}(U)\oplus H^{0,2}(U)) \cap 
H^2(U, {\bold R})$.

Since $q$ is positive-definite on the 
first two factors, $q$ has signature 
$(3, B-3)$. Therefore Corollary 8 has 
been proved 
when $X$ is projective.  

When $X$ is non-projective, by 
Proposition 7, we can 
find a point $s$ in any small open neighborhood 
$V$ of $0 \in S := \mathrm{Def}(X)$ 
in such a way that ${\cal X}_s$ is 
a projective symplectic variety. 
Since, for each $s \in S$, ${\cal X}_s$ 
is a symplectic 
variety with $h^0({\cal U}_s, \Omega^2_{{\cal U}_s}) 
= 1$, we can find 
a symplectic form $\Omega_s$ on each 
${\cal X}_s$ so that 
$\int_{{\cal U}_s}(\Omega_s {\overline \Omega_s})^l 
= 1$. Let 
$q_s$ be the 
quadratic form on $H^2({\cal U}_s, {\bold R})$ 
defined by $\Omega_s$. 
By using a natural flat structure in $R^2\pi_*{\bold C}
\otimes_{\bold C}{\cal O}_S$ ([Na, Theorem 8,(1)]), 
$H^2(U, {\bold R})$ and $H^2({\cal U}_s, {\bold R})$ 
are identified. By the same argument as 
[Be, Th{\'e}or{\`e}me 5], 
we see that $q$ and $q_s$ are proportional (by 
a positive constant) under this identification. 
Because ${\cal X}_s$ is projective, 
$q_s$ has signature $(3, B-3)$. Hence $q$ also 
has signature $(3, B-3)$. \vspace{0.15cm}  
 
{\bf Theorem 9}. {\em Notation and 
assumptions being the same as above, 
let $0 \in S_1 \subset S$ be a 
positive dimensional non-singular subvariety 
Then, for any open 
neighborhood $0 \in V \subset S$, 
there is a point $s \in V \cap S_1$ 
such that ${\cal X}_s$ is a projective 
symplectic variety.}  \vspace{0.15cm}

{\em Proof}. It suffices to prove 
that $N = 0$ in Proposition 7. 
But this follows from Corollary 8. 
\vspace{0.2cm}

\begin{center}
{\bf Appendix: Supplement to [Na]} 
\end{center}
  
In this appendix, we shall claim that 
Theorem 4 of [Na] remains true under 
a weaker condition. The exact statement 
is the following. \vspace{0.12cm}

{\bf Theorem A-4}. {\em  Let $X$ be a 
Stein open subset 
of a complex algebraic variety. 
Assume that $X$ has only rational singularities. 
Let $\Sigma$ be the singular locus of $X$ 
and let $f : Y \to X$ be a 
resolution of singularities 
such that $f\vert_{Y\setminus f^{-1}(\Sigma)} 
: Y \setminus f^{-1}(\Sigma) 
\cong X \setminus \Sigma$. 
Then $f_*\Omega^2_Y \cong i_*\Omega^2_U$ 
where $U := X \setminus \Sigma$ and 
$i : U \to X$ is a natural injection.} 
\vspace{0.15cm}

Theorem 4 in [Na] was stated under the condition 
that $X$ has rational Gorenstein singularities. 
We shall roughly sketch how to modify the 
original proof 
to drop the Gorenstein condition.  

The first step is to drop the Gorenstein 
condition from Proposition 1 of [Na]: 
\vspace{0.12cm} 

{\bf Proposition A-1}. {\em 
Let $X$ be a Stein open subset 
of a complex algebraic variety.  
Assume that $X$ has only rational 
singularities. 
Let $\Sigma$ be the singular locus of $X$ and 
let $f : Y \to X$ be a resolution of 
singularities 
such that 
$f\vert_{Y\setminus f^{-1}(\Sigma)} : 
Y \setminus f^{-1}(\Sigma) 
\cong X \setminus \Sigma$ 
and $D := f^{-1}(\Sigma)$ is a 
simple normal crossing 
divisor. 
Then $f_*\Omega^2_Y(\log D) 
\cong i_*\Omega^2_U$ 
where $U := X \setminus \Sigma$ and 
$i : U \to X$ 
is a natural injection.}  \vspace{0.15cm}

To prove Proposition 1, 
we have first taken an element 
$\omega$ from $H^0(U, \Omega^2_U)$, and 
have shown that $\omega$ has at worst a log 
pole at each irreducible component $F$ of $D$. 
When $X$ has rational Gorenstein singularities, 
$(X, p) \cong (R.D.P) \times 
({\bold C}^{n-2}, 0)$ 
for all singular points $p \in X$ 
outside certain codimension 3 (in $X$) locus 
$\Sigma_0 \subset \Sigma$. Since the proposition 
holds around such points, we only had to consider 
the case $f(F) \subset \Sigma_0$. 

In our general case, 
we have to take all irreducible 
components $F$ of $D$ into consideration. 
So we put 
$k := \dim \Sigma - \dim f(F)$, $l := 
\mathrm{Codim}(\Sigma \subset X)$ 
and continue the same 
argument as [Na, Proposition 1]. 
Here we note that a general hyperplane 
section $H$ 
of $X$ has again rational singularities.
In proving Claim of (a-2), 
we have used the following 
vanishing (the notation being 
the same as [Na]): 

$ R^i\pi_*\Omega^{l-2}_{Y_t}(\log D_t)(-D_t) = 
R^i\pi_*\Omega^{l-1}_{Y_t}(\log D_t)(-D_t) = 
R^i\pi_*\omega_{Y_t} = 0$ for  
$i \geq l -1$ and for $t \in \Delta$.  
   
Except the following cases, these vanishings 
follow from [St]: \vspace{0.12cm}

$l = 3$: $R^2\pi_*\Omega^1_{Y_t}(\log D_t)(-D_t)$, 
\vspace{0.12cm}

$l = 2$: $R^1\pi_*\Omega^1_{Y_t}(\log D_t)(-D_t)$, 
$R^2\pi_*{\cal O}_{Y_t}(-D_t)$, $R^1\pi_*{\cal O}
_{Y_t}(-D_t)$. \vspace{0.12cm}

For theses exceptional cases, we can prove the 
vanishing by combining the method in 
the proof of [Na-St, Theorem (1.1)] and 
the fact that a rational 
singularity is Du Bois [Kov].         

In proving Claim of (b-2) 
we also need similar vanishings; 
but they are already contained 
in the above cases. 
\vspace{0.15cm}

Next we shall generalize Lemma 2 of [Na] 
as follows: \vspace{0.12cm}

{\bf Lemma A-2}. 
{\em Let $p \in X$ be a Stein open 
neighborhood of a point $p$ of a 
complex algebraic 
variety. Assume that 
$X$ is a rational singularity of 
$dim X \geq 2$. 
Let $f : Y \to X$ be a resolution of 
singularities of 
$X$ such that $E := f^{-1}(p)$ is a 
simple normal crossing divisor. 
Then 
$H^0(Y, \Omega^i_Y) \to 
H^0(Y, \Omega^i_Y(\log E))$ 
are isomorphisms for $i = 1, 2$.} 
\vspace{0.15cm}

Lemma 2 of [Na] was stated under the 
condition 
$\dim X \geq 3$. In the proof of 
[Na, Lemma 2] 
we have first shown that 
$H^3_E(Y, {\bold C}) \to 
H^3(Y, {\bold C})$ and 
$H^2_E(Y, {\bold C})
 \to H^2(Y, {\bold C})$ 
are both injective. 
Even when $\dim X = 2$, these are true. 
First note that when $\dim X = 2$, 
$(X, p)$ is an isolated singularity.  
By taking the dual of the first map, we 
get the map $H^1_E(Y, {\bold C}) \to 
H^1(E, {\bold C})$. Since $X$ is a rational 
singularity, $H^1(E, {\bold C}) = 0$; thus 
the dual map is surjective. 

The injectivity of the second map 
follows from the next observation: 
$H^2(Y, {\bold C}) \cong H^1(Y, 
{\cal O}^*_Y)\otimes {\bold C}$, and 
$H^2_E(Y, {\bold C}) \cong 
\oplus {\bold C}[E_i]$ where 
$E_i$ are irreducible components 
of $E$. 

The rest of the proof of 
Lemma A-2 is the same as 
[Na, Lemma 2]. 
 
The following remark is 
quite similar to [Na, Remark 
below Lemma 2]: \vspace{0.12cm} 

{\bf Remark A}. {\em In 
Lemma A-2, the map 
$$H^0(E, \Omega^i_Y/\Omega^i_Y(\log E)(-E)) 
\to H^0(E, \Omega^i_Y(\log E)/
\Omega^i_Y(\log E)(-E))$$ 
is surjective for $i = 1, 2$} \vspace{0.12cm} 

Finally we shall prove \vspace{0.12cm} 

{\bf Proposition A-3}. {\em Let $X$ be a 
Stein open subset
 of a complex algebraic variety. 
Assume that $X$ has only rational 
singularities. 
Let $\Sigma$ be the singular locus of $X$ and 
let $f : Y \to X$ be a resolution of 
singularities 
such that $D := f^{-1}(\Sigma)$ is 
a simple normal 
crossing divisor and such that 
$f\vert_{Y \setminus D} : Y \setminus D 
\cong X \setminus \Sigma$. 
Then $f_*\Omega^2_Y 
\cong f_*\Omega^2_Y(\log D)$.}  
\vspace{0.2cm} 
 
This is a generalization of 
Proposition 3 of [Na], in which 
the same result was stated 
under the condition that $X$ has 
rational Gorenstein singularities.  
  
To prove Proposition 3, we have first 
taken an element 
$\omega$ from $H^0(X, f_*\Omega^2_Y
(\log D))$, and 
have shown that $\omega$ is regular 
along each irreducible component $F$ of $D$. 
When $X$ has rational Gorenstein singularities, 
$(X, p) \cong (R.D.P) \times 
({\bold C}^{n-2}, 0)$ 
for all singular points $p \in X$ 
outside certain codimension 3 (in $X$) locus 
$\Sigma_0 \subset \Sigma$. 
Since the proposition 
holds around such points, 
we only had to consider 
the case $f(F) \subset \Sigma_0$. 

In our general case, we have to 
take all irreducible 
components $F$ of $D$ into consideration. 
So we put 
$k := \dim \Sigma - \dim f(F)$, $l := 
\mathrm{Codim}(\Sigma \subset X)$ and 
continue the same 
argument as [Na, Proposition 3]. 
The rest of the argument is similar to  
[Na, Proposition 3]. In the original 
proof we have used [Na, Remark below 
Lemma 2], but now we shall use 
Remark A.  \vspace{0.2cm} 

{\bf Correction to [Na]}. 
In the introduction of [Na], 
a conjecture has been posed 
as a generalization of Bogomolov 
splitting theorem. This conjecture 
should be :  

{\em Let $Y$ be a smooth projective 
variety over {\bf C} with Kodaira 
dimension 0. Then there is a finite 
cover $\pi : Y' \to Y$ such 
that (a) $\pi$ is etale outside 
the support of the pluri-canonical 
divisor of $Y$, and (b) $Y'$ is 
birationally equivalent to 
$Y_1 \times Y_2 \times Y_3$, where 
$Y_1$ is an Abelian variety, $Y_2$ 
is a symplectic variety, and $Y_3$ 
is a Calabi-Yau variety.} 

In the old version, $\pi$ was 
assumed to be a finite etale 
cover. 

\begin{center}
Department of Mathematics, 
Graduate school of science, 
Osaka University, Toyonaka, Osaka 560, Japan 
\end{center}  
 
\end{document}